\documentclass[12pt]{article}
\usepackage{amsfonts,amsthm,amssymb}
\usepackage{graphics,color,url}
\usepackage{graphicx}
\usepackage{epsfig}

\newtheorem{thm}{Theorem}
\newtheorem{lem}{Lemma}

\newtheorem*{rem}{Remark}

\theoremstyle{definition}
\newtheorem{definition}{Definition}

\def\q#1.{{\bf #1.}}

\title{On proper colorings of hypergraphs}

\author{D.\,Karpov \and   N.\,Gravin}

\date{}

\begin{document}

\maketitle 

\thispagestyle{empty}

\section {Introduction}
In this paper we consider undirected graphs and hypergraphs. We denote by $V(G)$ the vertex set of a graph~$G$ 
and the edge set by~$E(G)$. Notations~$v(G)$ and~$e(G)$ in our paper stand for the number of vertices 
and edges respectively.

We denote by~$d_G(v)$ the {\it degree} of  vertex~$v\in V(G)$ in $G$. 
We denote by~$\delta(G)$ and~$\Delta(G)$ the minimal and maximal vertex degrees of~$G$ respectively.
We use similar notations ($V({\cal H})$, $E({\cal H})$ and~$d_{\cal H}(v)$) for a hypergraph~$\cal H$. 
In this work it is convenient for us to deal with edges and hyperedges in terms of vertex subsets 
of a graph or a hypergraph.

We denote the {\it neighborhood} of  vertex~$v$ in $G$ (i.e., the set of all adjacent to~$v$ vertices of $G$) by~$N_G(v)$. 

For any set $W\subset V(G)$ we denote by $G(W)$ the {\it induced subgraph} of~$G$ 
on~$W$ (i.e., the subgraph on~$W$ that contains all edges of~$G$ 
with two ends in~$W$).

There are several ways to generalize the notation of proper coloring on hypergraphs. For example, strong vertex colorings~\cite{AH}, in which all vertices in every hyperedge have to receive different colors. In the present 
paper we work with the definition proposed by P.\,Erd\H{o}s.

\begin{definition}
A vertex coloring of a hypergraph~$\mathcal{H}$ is called {\it proper coloring}, if any hyperedge contains at least two
vertices of different colors.
\end{definition}

In the field of colorings of ordinary graphs many natural questions are still left open. Thus, it is not surprising that
vertex colorings of hypergraphs are not well studied. 
One particular question in the field of hypergraph colorings that has received great attention in the literature (see~\cite{Be77,Be78,E63,E64,EL,Ko,Pl,Sc,Sp}) is the problem of finding a $n$-uniform hypergraph with the minimal number of edges~$m_k(n)$ that admits no proper vertex $k$-coloring. Another problem closely related to the one cited above is the question ``what is the minimal~$n$, such that every $n$-uniform and $n$-regular hypergraph (i.e., a hypergraph with all edges containing $n$ vertices and all vertices having degree~$n$) admits a proper vertex $2$-coloring?'' 
It was shown (see~\cite{AS}) by means of Lovasz local Lemma and other probabilistic methods, that for~$n\ge 9$ every such graph is $2$-colorable. Alon and Bregman~\cite{AB} improved this statement to~$n=8$, and Thomassen~\cite{Th} has shown finally 
$2$-colorablity for all~$n\ge 4$. 

The following theorem is the main result of our paper.

\begin{thm}
\label{th1}
Let $\mathcal{H}$ be a hypergraph of maximal vertex degree~$\Delta$, such that each its hyperedge
contains at least~$\delta$ vertices. Let $k=\lceil\frac{2\Delta}{\delta}\rceil$. 
Then the following statements hold.

$1)$ The hypergraph $\mathcal{H}$ admits proper vertex coloring in~$k+1$ colors.

$2)$ The hypergraph $\mathcal{H}$ admits proper vertex coloring in~$k$ colors, 
if~$\delta\ge 3$ and~$k\ge 3$.
\end{thm}

Our theorem gives weaker results than the works cited above, when the minimal size of hyperedge is close to the maximal vertex degree 
of considered hypergraph. However, for relatively small values of $\delta$ with respect to $\Delta$ the statement of our theorem 
becomes interesting. Our proof uses only classic combinatorial methods.

From our main theorem we derive results on dynamic vertex colorings.

\begin{definition}
A vertex coloring of a graph~$G$ is called {\it dynamic}, if any vertex~$v$ of degree at least~$2$ 
has at least two vertices of different colors in its neighborhood.
\end{definition}

We note that some papers (e.g.,~\cite{HMP,G,K}) study {\it proper} dynamic colorings. 
There it was shown the existence of a proper dynamic vertex coloring of~$G$ in 
$\Delta(G)+1$ colors~\cite{HMP} and in~$\Delta(G)$ colors~\cite{K} besides explicitly described 
series of exceptions. In the current paper we do not require a dynamic coloring to be a proper coloring
and obtain the following result. 

\begin{thm}
\label{cd}
Let~$G$ be a graph, $k=\lceil \frac{2\Delta(G)}{\delta(G)}  \rceil$. Then the following statements hold.

$1)$ The graph $G$  admits a dynamic vertex coloring in~$k+1$ colors.

$2)$ The graph $G$  admits a dynamic vertex coloring in~$k$ colors, if  $\delta(G)\ge 3$ and~$k\ge 3$. 
\end{thm}

\section{Hypergraph's Image and Alternating Chains}

We further introduce some notions which are important for the following proof of our main result.

\begin{definition}
We call any graph~$G$ (with possible multiple edges) by an {\it image} of a hypergraph~$\mathcal{H}$, 
if 

(i) $V(G)=V({\mathcal H})$;  

(ii) there exists a bijection $\varphi: E(G)\to E({\mathcal H})$, 
such that $e\subset \varphi(e)$ for every edge~$e\in E(G)$. 

We call~$\varphi$ by the {\it bijection of image}~$G$.
\end{definition}

\begin{rem}
We consider multiple edges of a graph-image~$G$ that corresponds to distinct hyperedges of the 
hypergraph~$\cal H$ as distinct edges.
\end{rem}

As in some classic theorems about vertex colorings, we make use of alternating chains. 
Next definition will show what we mean by this notion for hypergraphs.

\begin{figure}[!hb]
	\centering
		\includegraphics[width=0.8\columnwidth, keepaspectratio]{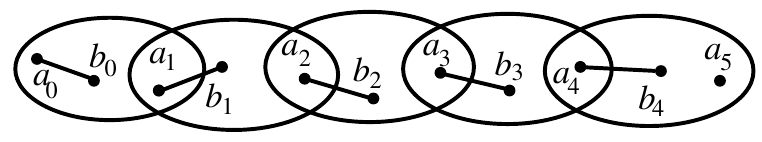}
     \caption{Alternating chain of length~5 with beginning~$a_0$ and end~$a_5$. }
	\label{fig1b}
\end{figure}

\begin{definition} Let $\delta \ge 3$ and let~$G$ be an image of a hypergraph~$\cal H$. We consider 
a sequence of vertices~$a_0b_0a_1b_1\dots a_n$ of~$\cal H$, satisfying the following conditions.
\begin{itemize}
\item For each $i$ vertices $a_i,b_i,a_{i+1}$ are different. 
\item There exist different  hyperedges~$e_0,\dots,e_{n-1}\in E({\cal H})$, such that
\begin{itemize}
\item $a_ib_i\in E(G)$ and $\varphi(a_ib_i)=e_i$,
\item $a_i,b_i,a_{i+1}\in e_i$. 
\end{itemize}

\end{itemize}

Then~$a_0b_0a_1b_1\dots a_n$ is an {\it alternating chain} from~$a_0$ to~$a_n$. We say that 
it has {\it length} $n$ and that it {\it goes} through the vertices~$a_0$, $b_0$,
\dots, $a_n$ and through the edges~$a_0b_0$,\dots, $a_{n-1}b_{n-1}$. We say that this chain 
{\it begins} at $a_0$ and {\it ends} at $a_n$.

For two sets $X,Y\subset V(G)$ with~$a_0\in X$ and~$a_n\in Y$, we say that
$a_0b_0a_1b_1\dots a_n$ is an alternating chain from~$X$ to~$Y$.
\end{definition}

\begin{rem}
$1)$  We allow the case~$n=0$ in the definition of alternating chain, that is~$a_0$ is an 
alternating chain from~$a_0$ to~$a_0$ of length~$0$.

$2)$ Since~$\varphi$ is a bijection, then edges~$a_1b_1,\dots, a_nb_n$ due to the definition of alternating chain
are all different. We recall here that multiple edges corresponding to different hyperedges of~$G$ are 
considered as different edges.

$3)$ Vertices are not necessarily different. An alternating chain may go through some vertices more than once.

\end{rem}

\begin{lem}
\label{obraz}
Let $\mathcal{H}$ be a hypergraph of maximal vertex degree~$\Delta$, such that each hyperedge of~$\mathcal{H}$
contains at least~$\delta$ vertices. Let $k=\lceil \frac{2\Delta}{\delta}\rceil$. Then there is an image~$G$
of~$\cal H$ with~$\Delta(G)\le k$.
\end{lem}

\renewcommand*{\proofname}{\bf Proof}
\begin{proof} Consider a trivial case~$\delta =2$. In this case for any image~$G$ of the hypergraph~$\cal H$ it is clear, that~$\Delta(G)\le \Delta=k$. In what follows we assume~$\delta \ge 3$.

For a graph~$H$ we denote by~$V_{k+1}(H)$ the set of all its vertices of degree at least~$k+1$.
We denote by~$s_{k+1}(H)$ the sum of degrees in the graph~$H$ taken over vertices of~$V_{k+1}(H)$.

For the sake of contradiction, we assume that the statement of lemma fails. Then for any image~$H$ we have~$V_{k+1}(H)\ne \varnothing$ 
and~$s_{k+1}(H)>0$.  Let~$G$ be an image with the minimal~$s_{k+1}(G)$. We denote by~$\varphi$ the bijection of~$G$, 
and we set $S=V_{k+1}(G)$.

Let~$U$ be the set of vertices of~$G$ that consists of all possible ends of alternating chains with the beginning in~$S$. 
We set~$F=G(U)$. Clearly, $U\supset S$. In the next we observe some properties of~$U$. 

\smallskip
\q1. {\it For any edge~$e\in E(F)$ the hyperedge~$\varphi(e)\subset U$.}

\noindent
Suppose the contrary. Then~$e=uw\in E(F)$ and the hyperedge~$\varphi(e)$ contains a vertex~$v\not\in U$
(see figure~\ref{fig2b}a). In the following we construct an alternating chain from~$S$ to~$v$ and, therefore, show that~$v\in U$.
The latter contradicts our assumption.

We consider the shortest alternating chain~$P=a_0b_0\dots a_n$ from~$S$ to~$\{u,w\}$. Without loss of generality 
we may assume that~$a_n=u$. Then $a_i\notin\{u,w\}$ for any $0\le i<n$. Hence~$P$ does not go through~$e=uw$. 
We add to~$P$ vertices~$w,v$ and obtain an alternating chain from~$a_0\in S$ to~$v$. 

\smallskip

\q2. {\it If a vertex~$u\in U$ is adjacent to a vertex~$v\notin U$, then all
the vertices of the hyperedge~$\varphi(uv)$ except~$v$ belong to~$U$.}

\noindent  Let~$u\in U$, $v\notin U$, $uv\in E(G)$, and~$e=\varphi(uv)$ be a hyperegde of~$\cal H$.
Suppose the contrary. We assume that $e$ contains a vertex~$w\notin U$ (see figure~\ref{fig2b}b).

As in the previous item, we construct the shortest alternating chain~$P$ from~$S$ to~$u$ 
(in the case~$u\in S$ this chain consists of one vertex). Let~$P$ goes through the edge~$uv$. 
Since we have chosen the shortest chain to~$u$, we have~$v=a_i$, $u=b_i$ for some $i<n$. 
Then we have~$v\in U$ and we arrive at a contradiction.

Thus the chain~$P$ does not go through the edge~$uv$. We add~$v$ and~$w$ to~$P$ and obtain
$w\in U$, that contradicts our assumption. Hence $v$ is the only vertex of the hyperedge~$e$ 
that does not belong to~$U$ (see figure~\ref{fig2b}c).

\begin{figure}[!hb]
	\centering
		\includegraphics[width=0.9\columnwidth, keepaspectratio]{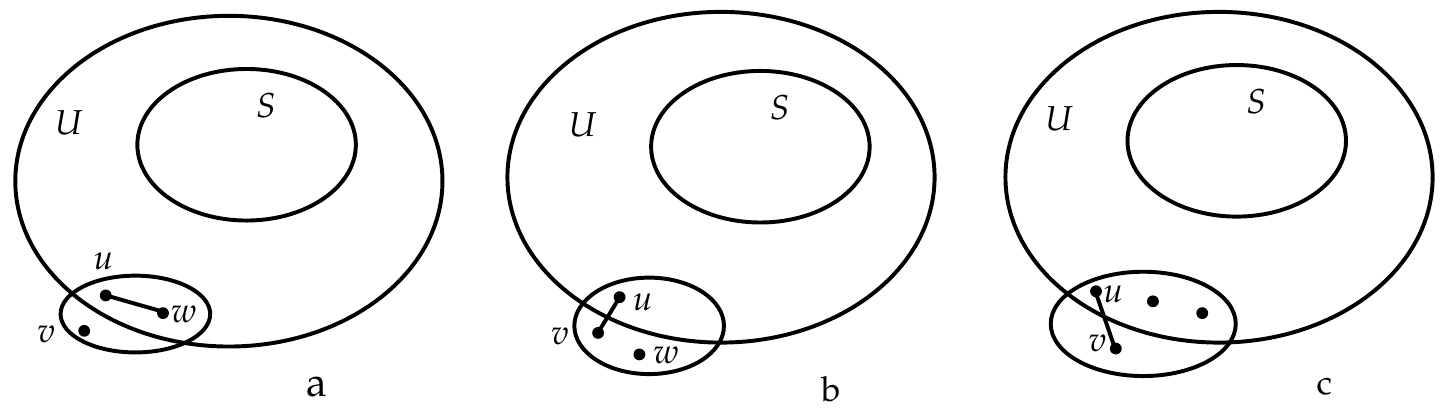}
     \caption{Hyperedges, intersecting~$U$. }
	\label{fig2b}
\end{figure}

\q3. {\it For any vertex~$u\in U$ we have~$d_G(u)\ge k$.}

\noindent Let~$u\in U$ and~$d_G(u)\le k-1$. Clearly,~$u\notin S$. 
Consider an alternating chain~$P=a_0b_0\dots a_n$ from~$S$ to~$u=a_n$.
We construct a new graph~$G'$: take the graph~$G$ and for every~$i\in [0,n-1]$ replace in the 
hyperedge~$e_i\supset \{a_i,b_i,a_{i+1}\}$ the edge~$a_ib_i$ by the edge~$b_ia_{i+1}$. 
It is easy to see that resulting graph~$G'$ is also an image of the hypergraph~$\cal H$.

Since~$d_{G'}(u)=d_{G}(u)+1\le k$, then~$u\not \in V_{k+1}(G')$. For any other vertex~$x$ we have
$d_{G'}(x)\le d_G(x)$. Hence $V_{k+1}(G')\subseteq S= V_{k+1}(G)$. It remains to notice that~$a_0\in S$ 
and~$d_G(a_0)> d_{G'}(a_0)$, consequently, $s_{k+1}(G')<s_{k+1}(G)$.
We obtain a contradiction with the minimality of~$s_{k+1}(G)$.

\smallskip
\q4. {\it We estimate the sum of degrees in the hypergraph~$\cal H$ taken over vertices in~$U$.}

\noindent 
Let~$u_1,\dots u_\ell$ be all vertices of~$U$ that have degrees less than~$k$  in 
the induced subgraph~$F=G(U)$. We set $$t_i=d_G(u_i)-d_F(u_i), \quad t=\sum_{i=1}^\ell t_i.$$ 
The degree of any vertex of~$U$ in~$G$ is at least~$k$ due to item~3. 
Since~$S\subset U$, the set~$U$ contains vertices that have degrees more than~$k$ in~$G$. 
Hence 
$$e(F) = {1\over 2}\sum_{u\in U} d_F(u) >{k|U| - t\over 2}. $$
We further estimate~$m=\sum_{u\in U} d_{\cal H}(u)$.
By item~1, all hyperedges of~$\mathcal H$ that correspond to the edges of~$F$ 
(i.e.  of the set $\varphi(E(F))$), are contained in the set~$U$ and contribute to~$m$ at least
$$\delta \cdot e(F) > \delta \cdot {k|U|  - t\over 2} \ge \Delta |U| - {\delta t \over 2}.$$

Now we consider~$t$ edges of~$G$ between~$U$ and~$V({\cal H})\setminus U$. According to item~2 each of these edges 
is contained in a hyperedge of~$\cal H$, which has only one vertex outside~$U$, and, consequently, 
which has at least~$\delta-1$ vertices in~$U$. We note that all these~$t$ hyperedges are different.
Thus
$$m> \Delta |U| - {\delta t \over 2} + (\delta-1) t> \Delta |U|.$$
Hence there is a vertex~$u\in U$ of degree~$d_{\mathcal{H}}(u)>\Delta$, that contradicts to the 
conditions of the lemma.

\smallskip
The obtained contradiction shows that there exists an image~$G$ of~$\cal H$ 
with~$\Delta(G)\le k$.
\end{proof}

\section{Proofs of Theorems~1 and~2}

\renewcommand*{\proofname}{\bf Proof of theorem~1}
\begin{proof} 
1)
We pick an image~$G$ of~$\cal H$ with~$\Delta(G)\le k=\lceil \frac{2\Delta}{\delta}\rceil$, which exists
due to the lemma~\ref{obraz}. Clearly, there exists a proper vertex coloring of the graph~$G$ in~$k+1$ 
colors. 

We need to show that this coloring is a proper vertex coloring of the hypergraph~$\cal H$. 
Let~$\varphi$ be the bijection of the image~$G$. For every hyperedge~${e\in E({\cal H})}$ we have~$\varphi^{-1}(e)\subset e$, 
and, therefore, two vertices of~$\varphi^{-1}(e) \subset E(G)$ have different colors.

2) To prove the second statement it suffices to find an image of the hypergraph~$\cal H$ that has a proper 
vertex coloring in~$k$ colors for~$k\ge 3$ and~$\delta\ge 3$.
At first we consider an image~$G$ of~$\cal H$ with~$\Delta(G)\le k$ and its bijection~$\varphi$. 

We remind the classic Brooks theorem:
{\it if $\Delta(G)\le k$,  $k\ge 3$, and no connected component of~$G$ is a clique on~${k+1}$ vertices, 
then~$G$ has a proper vertex coloring in~$k$ colors.} 

Let~$G$ have connected components that are cliques on~$k+1$ vertices. We enumerate them all by~$C_1$,\dots, $C_q$ 
(for conciseness, we will refer to these components simply by {\it cliques}).
Graph~$G$ can possibly have other connected components. We denote by~$D_{q+1},\dots, D_p$ induced subgraphs on 
 these components. In what follows we correct the graph-image~$G$, such that 
obtained graph would have proper vertex coloring in~$k$ colors. 

\smallskip
{\it Image transformation.}

\noindent Consider an arbitrary edge~$u_iw_i$ in each clique~$C_i$. It is clear that there is a
vertex~$v_i\in e_i=\varphi(u_iw_i)$ different from~$u_i$ and~$w_i$. We construct the new image~$G'$ of~$\cal H$, 
by replacing simultaneously every edge~$u_iw_i$ by the edge~$u_iv_i$. 
We call the edges~$u_1v_1$, \dots, $u_qv_q$ by {\it new edges}.

\smallskip
{\it Further we prove that~$G'$ has a proper vertex coloring in~$k$ colors.}

\noindent
We construct an auxiliary digraph~$F$: vertices of~$F$ are connected components of~$G$, 
from each component-clique~$C_i$ an oriented edge ({\it arc}) leads to a component that contains~$v_i$. 
If~$v_i$ is a vertex of the clique~$C_i$, then this arc will be a loop.
In fact, in order to construct $F$ from $G'$, one could orient the new edges and contract each connected 
component of~$G$ into a vertex.  

Our algorithm for coloring vertices in~$k$ colors works according to the following plan:

\smallskip
--- if there exists a clique that has no incoming arc in~$F$, we perform Step~1 
and return to the beginning of the algorithm;

\smallskip
--- if each clique has at least one incoming arc in~$F$, then we perform Step~2 and  
terminate the algorithm.

\smallskip
\q1. {\it There is clique~$C_i$ that has no incoming arc.}

\noindent In this case $d_{G'}(w_i)=k-1$. We enumerate vertices of~$C_i$ starting from~$w_i$ 
and finishing at the vertex~$u_i$ that is adjacent in~$G'$ to a vertex of another connected 
component of~$G$. We assume that vertices of the rest components are properly colored in~$k$ colors. 
Then we can color vertices of~$C_i$ in the reverse order (respect to their numbers): 
at each step we take a vertex that is adjacent to less than $k$ already colored vertices and 
we color it in any remaining color. 

Therefore, we can delete from~$G'$ all vertices of the component~$C_i$ and continue by
coloring the remaining graph ${G'- C_i}$. In addition, we change the graph~$F$. 
We delete from $F$ vertex~$C_i$ and the arc going from~$C_i$.

\smallskip
\q2. {\it Every component-clique has  an incoming arc.}

\noindent Since exactly one arc goes from each clique, then exactly one arc comes into each clique. 
Thus all cliques in $F$ are divided into several oriented cycles, which vertices are not adjacent 
to each other in~$G'$. We color these cycles independently. 
The rest connected components of~$G$ (not cliques on~${k+1}$ vertices) are the same connected components 
in~$G'$. Due to Brooks theorem their vertices can be properly colored in~$k$ colors.

Now we have cliques~$C_1$, \dots, $C_\ell$ forming in $F$ an oriented cycle. We denote by~$G^*$ 
the induced subgraph of~$G'$ on the union of all these cliques.  It remains to prove that~$G^*$ 
has a proper vertex coloring in~$k$ colors. If~$\Delta(G^*)\le k$, then it follows from Brooks theorem, 
since the graph~$G^*$ is connected and is not a clique on~${k+1}$ vertices.
Assume, that $\Delta(G^*)>k$ and consider two cases.

\smallskip
\q{2.1}. {\it $\ell=1$, i.e. our cycle is a loop and~$v_1\in V(C_1)$.}

\noindent Then $G^*$ is a clique on~${k+1}$ vertices with deleted edge~$u_1w_1$ and edge $u_1v_1$ of multiplicity two. Clearly,~$G^*$ has a proper coloring in~$k$ colors: we color~$u_1$ and $w_1$ in the same color, and we color each other 
vertex in its own color.

\smallskip
\q{2.2}. $\ell\ge 2$.

\noindent
Let~$G^*$ have a vertex~$x$ of degree more than~$k$ and $x$ belongs to the clique~$C_i$. 
Clearly, $x$ is adjacent to a vertex of the clique~$C_{i-1}$ and~$x\ne w_i$. Moreover, 
in this case~$d_G(w_i)=k-1$. We can delete the vertex~$w_i$ from~$G^*$, since we can color 
this vertex after coloring the rest vertices. If~$x\ne u_i$, then~$x$ is adjacent to~$w_i$, 
hence, all remaining vertices of~$C_i$ have degrees not more than~$k$ in~${G^*-w_i}$ 
(see figure~\ref{fig3b}a).

\begin{figure}[!hb]
	\centering
		\includegraphics[width=0.9\columnwidth, keepaspectratio]{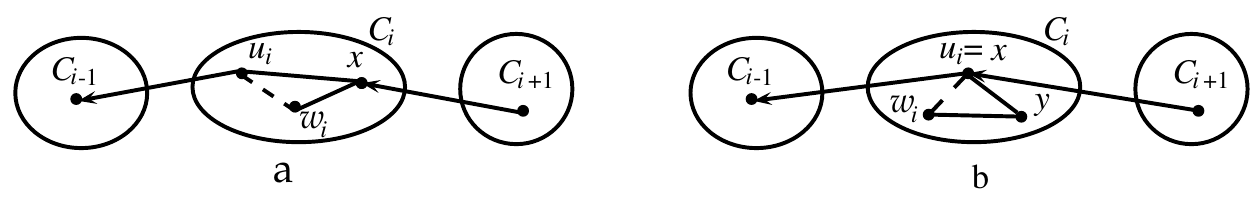}
     \caption{Coloring of the clique~$C_i$ in the graph~$G^*$.}
	\label{fig3b}
\end{figure}

If~$x=u_i$, then there is another vertex~$y$ in~$C_i-w$ and~$d_{G^*-w_i}(y)= k-1$. 
We delete~$y$ from the graph (after coloring all other vertices in $k$ colors we can easily color~$y$).
Clearly, in the graph~${G^*-w_i-y}$ the degrees of all remaining vertices of~$C_i$ do not exceed~$k$ (see figure~\ref{fig3b}b). 

We perform such operations with all components-cliques that have in~$G^*$ a vertex of degree more than~$k$. 
As a result we obtain a connected graph~$H^*$ with maximal degree not exceeding~$k$. 

\smallskip \goodbreak
{\it We prove that~$H^*$ is not a clique on~${k+1}$ vertices.}

\noindent It is clear from the construction that all new edges between the components~$C_1,\dots, C_n$ remain in~$H^*$ after deletion
described above (we have not deleted none of their ends from~$G^*$). We consider a component~$C_2$ and two new edges~$u_1v_1$ and~$u_2v_2$, incident to  vertices of~$C_2$. Clearly, the graph ${H^*-u_1v_1-u_2v_2}$ is disconnected  ($C_2$ is separated from the rest vertices of the graph). Thus $H^*$ becomes disconnected after deleting two of its edges, and, hence, it can't be a clique on~$k+1\ge 4$ vertices.

\smallskip
Applying Brooks theorem we get a proper $k$-coloring of~$H^*$.
After that we can add back all deleted from~$G^*$ vertices and color them properly 
in reverse order to their deletion.

Thus the graph $G'$ admits a proper $k$-coloring and this coloring, as it was mentioned above,
provides a proper vertex coloring of the hypergraph~$\cal H$.
\end{proof}

\renewcommand*{\proofname}{\bf Proof of theorem~2}
\begin{proof}

We construct the following hypergraph~$\mathcal{H}$. 
Its vertex set~$V(\mathcal{H})$ coincides with~$V(G)$; set of hyperedges~$E(\mathcal{H})$ consists of 
neighborhoods $N_G(v)$ of all vertices~$v\in V(G)$.
Each hyperedge of~$\mathcal{H}$ has the size at least~$\delta(G)$ and each vertex of~$\mathcal{H}$ belongs 
to not more than~$\Delta(G)$ hyperedges. Now it is easy to see, that the statement we are proving is an 
immediate consequence of theorem~\ref{th1} applied to the hypergraph $\mathcal{H}$.
\end{proof}

\newcounter{roman}
\setcounter{roman}{2}


\begin{thebibliography}{99}


\bibitem{AH} {\sc G.\,Agnarsson,  M.\,M.\,Halldo\'rsson.} 
{\it Strong Colorings of Hypergraphs.} Approximation and Online Algorithms.
Lecture Notes in Computer Science, 2005, Volume 3351/2005, p.253-266.

\bibitem{AB} {\sc N.\,Alon and Z.\,Bregman.} {\it Every $8$-uniform $8$-regular hypergraph is $2$-colorable.} 
Graphs Combinat. 4 (1988), p.303–305.

\bibitem{AS} {\sc N.\,Alon and J.\,Spencer.} {\it The probabilistic method.}
 Wiley-Interscience, New York, 2000.

\bibitem{Be77} {\sc J.\,Beck.} {\it On a combinatorial problem of P.\,Erd\H{o}s and L.\,Lov\'{a}sz.} Discrete Math. 17 (1977), p.127-131.

\bibitem{Be78} {\sc J.\,Beck.}  {\it On $3$-chromatic hypergraphs.}, Discrete Math. 24 (1978), p.127-137.

\bibitem{E63} {\sc P.\,Erd\H{o}s.} {\it On a combinatorial problem.}  Nordisk. Mat. Tidskr. 11 (1963), p.5-10.

\bibitem{E64} {\sc P.\,Erd\H{o}s.} {\it On a combinatorial problem, \Roman{roman}.} 
Acta Math. Acad. Sci. Hungar. 15 (1964), p.445-447.

\bibitem{EL} {\sc P.\,Erd\H{o}s,  L.\,Lov\'{a}sz.} 
{\it Problems and results on $3$-chromatic hypergraphs and some related questions.} Infinite and finite sets, Colloq. Math. Soc. J. Bolyai, Vol. 10, North Holland, Amsterdam, 1974, p.609-627.

\bibitem{HMP} {\sc L.\,Hong-Jian, B.\,Montgomery, H.\,Poon.} \textit{Upper
Bounds of Dynamic Chromatic Number.}  Ars. Combinatoria 68(2003), p.193-201.


\bibitem{Ko} {\sc A.\,Kostochka.} {\it Coloring uniform hypergraphs with few colors.} Random Structures and Algorithms 24 (2004), p.1-10.


\bibitem{Pl} {\sc A.\,Pluha'r.} {\it Greedy colorings of uniform hypergraphs.} Random Structures and Algorithms 35, (2009)  p.216–221.

\bibitem{Sc} {\sc W.\,M.\,Schmidt.} {\it Ein kombinatoriches problem/} Acta Math. Acad. Sci. Hungar 15 (1964), p.373-374.

\bibitem{Sp} {\sc J.\,H.\,Spencer.} {\it Coloring $n$-sets red and blue.} J. Combin Theory Ser. A 30 (1981), p.112-113.

\bibitem{Th} {\sc C.\,Thomassen.} {\it The even cycle problem for directed graphs.} J. Am.Math. Soc. 5 (1992), 
p.217-229.


\bibitem{G} {\sc N.\,V.\,Gravin.} {\it Nondegenerate colorings in the Brooks theorem.} 
Diskretn. Mat., {\bf 21} (2009), i.4, p.106-128, in Russian.
English   translation in Discrete Math. Appl.  {\bf 19}  (2009),  no. 5, p.533-553.  


\bibitem{K} {\sc D.\,V.\,Karpov.} {\it Dynamic proper vertex colorings of a graph.} 
 Zap. Nauchn. Semin. POMI v.381 (2010), p.47-77. English translation to appear in Journal of mathematical Sciences.



\end{thebibliography}
\end{document}